\documentclass[a4paper,12pt]{article}
\usepackage[utf8x]{inputenc}

\usepackage[normalsize,sc]{caption}
 \usepackage{epsf,epsfig,latexsym,amssymb}
 \usepackage{amsthm}
 \usepackage{cite,graphicx}
 \usepackage{multirow}

\usepackage{mathrsfs}
\usepackage{amsmath}

\theoremstyle{definition}
{\rm}
{\rm}

\usepackage[version=3]{mhchem} 

\usepackage[colorlinks=true, linkcolor=blue, citecolor=blue]{hyperref}
 
   {\begin{trivlist}\item[]{\bfseries\sffamily Keywords:}\ }
   {\end{trivlist}}
\usepackage{auto-pst-pdf}
\usepackage{pst-plot}
\usepackage{pdftricks}
\begin{psinputs}
\usepackage{pstricks,pst-plot}
\end{psinputs}

\usepackage{fancyhdr}

\begin{document}



\title{Optimal control of crystallization of alpha-lactose monohydrate\footnote{Draft of accepted article in IEEE {\it Asian Control Conference (ASCC 2013)}. The use of this material is subject to IEEE Copyright policy - 
\url{http://www.ieee.org/publications_standards/publications/rights/reqperm.html}}}
\author{A. Rachah\thanks{Universit\'e Paul Sabatier, Institut de Math\'ematiques, Toulouse, France} 
\and D. Noll$\!\,^\dag$
\and  F.  Espitalier\thanks{Ecole des Mines d'Albi, Centre Rapsodee, Albi, France}\and F. Baillon$\!\, \dag$
}
\date{}
\maketitle

\begin{abstract}
We present a mathematical model for solvated crystallization of $\alpha$-lactose 
monohydrate in semi-batch mode. The process dynamics are governed by conservation laws 
including population, molar 
and energy balance equations. We present and discuss the model and then 
control the process with the goal to privilege the production of small particles in the range 
between $10^{-5}$ and $10^{-4}\mu m$. 
We compare several specific and unspecific cost functions leading to
optimal strategies with significantly different effects on product quality. Control 
inputs are
 temperature, feed rate, and the choice of an appropriate
crystal seed.
\end{abstract}

\section{Introduction}
\label{intro}

Crystallization is the unitary operation of formation of solid crystals from a 
solution. In process engineering crystallization is an important separation 
process used in chemical, pharmaceutical, food,  material and semiconductor industries. 
Mathematical models are described by conservation laws with population, molar and energy
balance equations.
Crystallizers can be operated either in batch, semi-batch or continuous mode. Semi-batch 
crystallization is widely used in the pharmaceutical and fine chemical industry for the production
of solids in a variety of operating modes.

Crystallization processes are described by balance equations, including
a population balance for the particle size distribution, a molar balance
for the distribution of solute, and an energy balance equation to
model thermodynamic phenomena. 
In the food-processing industry, there has been a growing interest in the crystallization of lactose in recent years
 \cite{jer, ami, amishu}.
In this paper we study a model of solvated crystallization of $\alpha$-lactose monohydrate,
which includes four interacting populations, one of them aging, in tandem with an energy balance.
Two forms of lactose ($\alpha$- and $\beta$-lactose) exist simultaneously in
aqueous solution, the exchange being governed by mutarotation.

\hspace*{-1.2cm}

\begin{picture}(100,100)

\unitlength.042cm

\put(0,0){\line(0,1){30}}
\put(0,0){\line(1,0){50}}
\put(0,30){\line(1,0){50}}
\put(50,0){\line(0,1){30}}
\put(5,19){$\beta$-lactose}

\put(100,0){\line(0,1){30}}
\put(100,0){\line(1,0){50}}
\put(100,30){\line(1,0){50}}
\put(150,0){\line(0,1){30}}
\put(106,19){$\alpha$-lactose}

\put(68,26){$k_2$}
\put(68,-4){$k_1$}
\put(50,22){\vector(1,0){50}}
\put(100,8){\vector(-1,0){50}}



\put(200,-60){\line(0,1){30}}
\put(200,-60){\line(1,0){50}}
\put(200,-30){\line(1,0){50}}
\put(250,-60){\line(0,1){30}}
\put(211,-49){water}

\put(250,-50){\vector(1,0){30}}
\put(280,-40){\vector(-1,0){30}}
\put(282,-53){evap.}
\put(282,-43){feed}

\put(225,-30){\vector(0,1){30}}
\put(231,-19){hydration}
\put(150,15){\vector(1,0){50}}
\put(162,20){nucl.}

\put(200,0){\line(1,0){50}}
\put(200,0){\line(0,1){30}}
\put(250,0){\line(0,1){30}}
\put(200,30){\line(1,0){50}}
\put(208,18){$\alpha$-mono}
\put(208,5){crystal}
\put(280,15){\vector(-1,0){30}}
\put(284,13){seed}

\put(200,60){\line(1,0){50}}
\put(200,60){\line(0,1){30}}
\put(250,60){\line(0,1){30}}
\put(200,90){\line(1,0){50}}
\put(208,78){$\alpha$-mono}
\put(208,65){crystal}
\put(250,80){\vector(1,0){30}}
\put(282,78){product}
\put(280,70){\vector(-1,0){30}}
\put(284,66){seed}

\put(235,30){\vector(0,1){30}}
\put(215,60){\vector(0,-1){30}}
\put(195,43){attr.}
\put(239,43){growth}

\put(25,60){\vector(0,-1){30}}
\put(0,55){feed}

\put(125,60){\vector(0,-1){30}}
\put(100,55){feed}


\end{picture}
\vspace*{2.2cm}
\begin{figure}[h!]
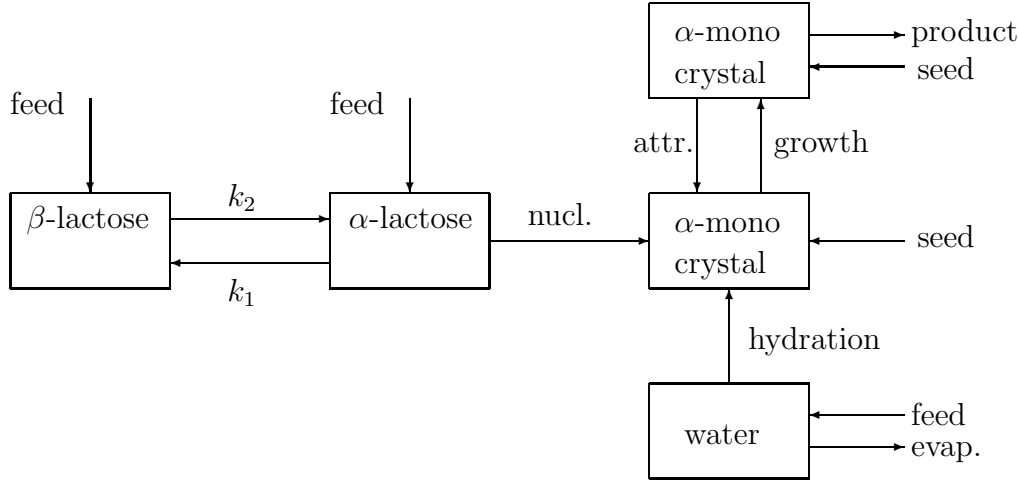

\caption{Solvated crystallization of $\alpha$-lactose monohydrate with complex population dynamics featuring
one aging and three ageless populations. Exchange rates depend on temperature. Controls are
feed rate in semi-batch mode, temperature of envelope, and the distribution of the crystal seed.
Based on mathematical modeling the process is optimized
in order to maximize the particle mass in a given small
size range  $L_{\rm low}\leq L \leq L_{\rm high}$. \label{fig-lactose}}
\end{figure}

Crystallization and precipitation processes are modeled as highly nonlinear and complex dynamical
systems. This makes it interesting
to simulate,
control and optimize  these processes in order  to enhance product quality in various situations
 \cite{rand,james}. 
 In order to  control crystallization processes, Hu and Rohani \cite{hu} have studied different
 heuristic cooling methods such as
 linear cooling, natural cooling and controlled cooling.  
In this study we   control of the process of formation of $\alpha$-lactose
crystals in such a way that the growth of very large crystals is avoided and the bulk of
crystal mass occurs in a small particle range between $10^{-5}$ and $10^{-4}\mu m$. 
In order to achieve this goal we use optimal control of the process in semi-batch mode by variations of
temperature, feed rate, and also by an appropriate choice of the crystal seed, and based on
a variety of different figure of
merit functions.

\section{Dynamic model of process}
\label{dynamic}
In this section
the dynamic model of semi-batch crystallization of $\alpha$-lactose
monohydrate  is described by presenting the  population, molar and energy
balance equations. 

\vspace{.1cm}\noindent
{\bf Population balance equation:}
The population balance equation is a first-order PDE
\begin{eqnarray}
\label{pop}
\frac{\partial}{\partial t} \left( V(t) n(L,t) \right) + V(t) G\left(c_\alpha(t),c_\beta(t),T(t)\right) \frac{\partial n(L,t)}{\partial L}
= V(t)  \dot{n}(L,t)^\pm,
\end{eqnarray}
where  $n(L,t)$ is the distribution of $\alpha$-lactose crystals, 
$c_\alpha(t),c_\beta(t)$ are the dimensionless
concentrations of $\alpha$- and $\beta$-lactose in the liquid phase,
$V(t)$ is the volume of slurry in the crystallizer, a dependent
variable given in (\ref{V}), $G\left(c_\alpha,c_\beta,T\right)$ is the temperature-dependent 
growth coefficient of $\alpha$-crystals, assumed independent of crystal size $L$, and the
right hand side represents source and sink terms.  We add the boundary condition
\begin{eqnarray}
\label{bdry}
n(0,t) = \frac{B\left(c_\alpha(t),c_\beta(t),T(t)\right)}{G\left( c_\alpha(t),c_\beta(t),T(t)\right)},
\end{eqnarray}
and the initial condition
\begin{eqnarray}
n(L,0) = n_0(L),
\end{eqnarray}
where $n_0(L)$ is the crystal seed.
It is convenient to introduce the
moments of the crystal size distribution function
\[
\mu_\nu(t) = \int_0^\infty n(L,t)L^\nu dL, \quad \nu = 0,1,\dots,
\]
which allows to break (\ref{pop}) into an infinite sequence of ODEs
if the source and sink terms may be neglected. In the present situation this amounts to neglecting
agglomeration and also attrition effects.
Using (\ref{bdry}) we obtain the equations
\begin{eqnarray}
\label{mom1}
\frac{d\mu_\nu(t)}{dt} + \frac{V'(t)}{V(t)} \mu_\nu(t) - \nu G\left( c_\alpha(t),c_\beta(t),T(t)\right) \mu_{\nu-1}(t) = 0,
\end{eqnarray}
$\nu=1,2,\dots$, in tandem with
\begin{eqnarray}
\label{mom2}
\frac{d\mu_0(t)}{dt} + \frac{V'(t)}{V(t)} \mu_0(t) - B\left( c_\alpha(t),c_\beta(t),T(t)\right)=0.
\end{eqnarray}
The initial conditions are then
\[
\mu_\nu(0) = \int_0^\infty n_0(L)L^\nu dL, \quad \nu=0,1,\dots. 
\]

\vspace{.1cm}\noindent
\textbf{Solvent mass balance:}
\begin{eqnarray}
\label{water}
\frac{dm_{\rm H_2O}(t)}{dt}& = (R^{-1}-1) 3 k_v \rho_{\rm cry} G\left(c_\alpha(t),c_\beta(t),T(t)  \right)  
 V(t) \mu_2(t) + q_{\rm H_2O}(t)
\end{eqnarray}
Here $m_{\rm H_2O}$ is the mass of water in the aqueous solution, which
changes due to feed, $q_{\rm H_2O}$, and due to the integration of water molecules into the
$\alpha$-crystals. The constant $R=M_{\rm cry}/M_{\alpha}=1.0525$ is the ratio of the molar masses
of the solid and liquid phases of $\alpha$-lactose.

\vspace{.1cm}\noindent
\textbf{Concentration of $\alpha$-lactose:} The dimensionless concentration
of $\alpha$-lactose $c_\alpha$  in the solution is defined as $m_\alpha =c_\alpha m_{\rm H_2O}$ and satisfies
the differential equation
\begin{eqnarray}
\label{alpha}
\frac{dc_\alpha(t)}{dt}=
& 
\displaystyle
\frac{1}{m_{\rm H_2O}(t)}[c_\alpha(t) (1-R^{-1})-R^{-1}] \frac{dm_{\rm cry}(t)}{dt}  
-k_1(T(t))c_\alpha(t)+k_2(T(t))c_\beta(t)\\
\displaystyle
&\notag + \left( c_\alpha^+(t)-c_\alpha(t)\right) \displaystyle \frac{q_{\rm H_2O}(t)}{m_{\rm H_2O}(t)}.
\end{eqnarray}
Here $m_{\rm cry}$ is the crystal mass in the slurry, $c_\alpha^+$ is the feed rate of $\alpha$
lactose during the semi-batch phase, and $k_1(T)$, $k_2(T)$ are the temperature dependent
mutarotation exchange rates between $\alpha$- and $\beta$-lactose in the liquid phase.

\vspace{.1cm}\noindent
\textbf{Concentration of $\beta$-lactose:} The dimensionless concentration of $\beta$-lactose $c_\beta$ 
is defined as $m_\beta = c_\beta m_{\rm H_2O}$ and satisfies the differential equation
\begin{eqnarray}
\label{beta}
\frac{dc_\beta(t)}{dt}&=\displaystyle
\frac{c_\beta(t)}{m_{\rm H_2O}(t)}(1-R^{-1}) \frac{dm_{\rm cry}(t)}{dt} + k_1(T(t))c_\alpha(t) - k_2(T(t))c_\beta(t)  \\
&\notag +\left( c_\beta^+(t)-c_\beta(t)\right)\displaystyle \frac{q_{\rm H_2O}(t)}{m_{\rm H_2O}(t)}.
\end{eqnarray}

\vspace{.1cm}\noindent
\textbf{Energy balance:}
The temperature hold system describes the interaction between
crystallizer temperature, the temperature of the jacket, and the control signal,
to which the internal heat balance due to enthalpy is added. We have
\begin{align}\label{eq3}
\dfrac{dT(t)}{dt}&=P_1(t)\bigg[ -P_2(t)(T(t)-T_{\text{ref}})-\Delta H \frac{dm_{\text{cry}}(t)}{dt}
+ U A(t) \left(T_{\text{jacket}}(t)-T(t)\right)\\
&\!\!\!\!\!\!\!\!+q_{{\rm H_2O}}(t)\left(C^p_{{\rm H_2O}}+C^p_{\alpha}c_{\alpha}(0)+C^p_{\beta}c_{\beta}(0)\right)\left(T_{\text{feed}}-T_{\text{ref}} \right)\bigg]
\notag
\end{align}
where
\begin{eqnarray}\label{eq4}
\dfrac{dT_{\text{jacket}}(t)}{dt}=-0.0019(T_{\text{jacket}}(t)-T_{\rm sp}(t))
\end{eqnarray}
was obtained through identification of the system.
Here $T(t)$ is the temperature of the slurry, $T_{\text{ref}} =25^0$C a constant
reference temperature, $T_{\text{feed}}$ the temperature of feed, 
which is the temperature of ${\rm H_2O}$ in this case, assumed constant in this study, 
$T_{\rm jacket}(t)$ is the temperature of the
crystallizer jacket, and $T_{\rm sp}(t)$ is the   set point temperature, which is used
as a control input to regulate $T_{\rm jacket}(t)$, and therefore indirectly
$T(t)$ via the heat exchange between the envelope
and the crystallizer through the contact surface $A(t)$, which is determined through
$V(t)$.  The constants $C^p_{\rm H_2O}$, $C_\alpha^p$, $C_\beta^p$ 
are the specific heat capacities.
We have
\begin{align*}
P_1(t)^{-1}=m_{\rm H_2O}(t)C^p_{\rm H_2O}+m_{\alpha}(t)C^p_{\alpha}+m_{\beta}(t)C^p_{\beta} +m_{\rm cry}(t)C^p_{\rm cry},
\end{align*}
\begin{align*}
P_2(t)=\dfrac{dm_{\rm H_2O}(t)}{dt}C^p_{\rm H_2O}+\frac{dm_\alpha(t)}{dt}C^p_{\alpha}+
\frac{dm_\beta(t)}{dt}C^p_{\beta} +\dfrac{dm_{\rm cry}(t)}{dt}C^p_{\rm H_2O}, 
\end{align*}
with $m_\alpha = c_\alpha m_{\rm H_2O}$, $m_\beta = c_\beta m_{\rm H_2O}$.


\vspace{.1cm}\noindent
\textbf{Mutarotation:}
The mutarotation exchange coefficients $k_1,k_2$ are temperature dependent and
are determined as

$$ k_2(T) = k_0\cdot \exp(-E/(R\cdot (T+273.15))),$$
 
$$ k_m(T) = 1.64-0.0027\cdot T, k_1(T) = k_2(T)\cdot k_m(T).$$

\noindent
The equilibrium of mutarotation therefore occurs at

$$c_{\alpha,\rm sat, eq}(T) = \dfrac{10.9109\cdot \exp(0.02804\cdot T)}{100(1+k_m(T))},$$

$$ F(T) = 0.0187\cdot \exp(0.0236\cdot T),$$
  
\begin{eqnarray*} c_{\alpha,\rm sat}(c_\beta,T) = c_{\alpha,\rm sat, eq}(T)-F(T) (c_\beta-k_m(T)\qquad \\\times c_{\alpha,\rm sat, eq}(T)).
\end{eqnarray*}

\vspace{.1cm}\noindent
\textbf{Nucleation rate:} The nucleation rate is based on a phenomenological law
\begin{eqnarray*}
B(c_\alpha,c_\beta,T)
=
k_b \exp \left(  - \frac{B_0}{(T+273.15)^3\ln^2\left( \frac{c_\alpha}{c_{\alpha,\rm sat}(c_\beta,T)}\right)}\right)
\end{eqnarray*}
as is the growth rate

\vspace{.1cm}\noindent
\textbf{Growth rate:}
\begin{eqnarray*}
G(c_\alpha,c_\beta,T)
=
k_g \left( c_\alpha - c_{\alpha,\rm sat}(c_\beta,T)\right).
\end{eqnarray*}

\vspace{.1cm}\noindent
\textbf{Volume:} The total volume of slurry $V(t)$ is a dependent variable, which can be expressed as a function of
the states $c_\alpha$, $c_\beta$ and $m_{\rm H_2O}$ through
\begin{align}
\label{V}
V(t) = \frac{m_{\rm H_2O}(t)}{1-k_v\mu_3(t)} \left[ \rho_{\rm lac,\alpha}^{-1} c_\alpha(t) + \rho_{\rm lac,\beta}^{-1} c_\beta(t) + \rho_{\rm H_2O}^{-1}\right].
\end{align}
Therefore
\begin{align*}
\dfrac{dV(t)}{dt}&= 3k_v G(c_\alpha(t),c_\beta(t),T(t)) V(t) \mu_2(t)
+\frac{dm_{\rm H_2O}(t)}{dt} \left[ \rho_{\rm lac,\alpha}^{-1} c_\alpha(t)+\rho_{\rm lac,\beta}^{-1} c_\beta(t)
 +\rho_{\rm H_2O}^{-1} \right] \\
& +m_{\rm H_2O}(t) \left(\rho_{\rm lac,\alpha}^{-1} \frac{dc_\alpha(t)}{dt} +\rho_{\rm lac,\beta}^{-1} \frac{dc_\beta(t)}{dt} \right).
\end{align*}

\vspace{.1cm}\noindent
\textbf{Crystal mass:} The total crystal mass satisfies the equation
\begin{eqnarray}
\label{dmcry}
\frac{dm_{\rm cry}(t)}{dt} =
3 k_v\rho_{\rm cry} G(c_\alpha(t),c_\beta(t),T(t)) V(t) \mu_2(t).
\end{eqnarray}


\begin{table}[h!]
\begin{center}
\hspace*{.10cm}
\begin{tabular}{||c|c|c|c||}
\hline\hline
$R$ & 1.0525 & -- &ratio of molar masses  \\
\hline
$k_v$ &  0.523598  & --   &volumic shape factor \\
\hline
$\rho_{\rm cry}$ & 1545 &   $kg \cdot m^{-3}$ &crystal density \\
\hline
$\rho_{\rm lac,\alpha}$& 1545 & $kg \cdot m^{-3}$  &$\alpha$-lactose density\\
\hline
$\rho_{\rm lac,\beta}$& 1590 & $kg \cdot m^{-3}$  &$\beta$-lactose density\\
\hline
$\rho_{\rm H_2O}$ &1000 & $kg \cdot m^{-3}$ & water density\\
\hline
$\Delta H$ & -43.1 & $kJ/kg$ &   heat of crystallization \\
\hline
$U$ &300&  $kJ/m^{2}.h.K$  & heat transfert coefficient  \\
\hline
$k_0$ & $2.25\cdot 10^{8}$ &$s^{-1}$&kinetic mutarotation\\ &&& constant  \\
\hline
$T_{\text{ref}}$ & 25 & ${\,}^0C$& reference temperature \\
\hline
$C^p_{\rm H_2O}$ &  4180.5&  $J/kg/K$ & heat capacity ${\rm H_2O}$\\
\hline
$C^p_{\rm cry}$ &1251 & $J/kg/K$ &heat capacity \\ &&&$\alpha$-lactose monohydrate\\
\hline
$C^p_\alpha$ &1193& $J/kg/K$ &heat capacity $\alpha$-lactose\\
\hline
$C^p_\beta$ &1193&$J/kg/K$ &heat capacity $\beta$-lactose\\
\hline
$B_0$ & 5.83 & & nucleation constant\\
\hline
$R_g$ & $18.314 $  & $J/K/mol$   &universal gas constant\\
\hline
$k_g$ & $10 \cdot 10^{10}$  & $m\cdot s^{-1}$  &growth rate coefficient\\
\hline
$k_b$ & $1.18 \cdot 10^{-7}$ &$\sharp \cdot m^{-3}s^{-1}$ &birth rate coefficient\\
\hline
$t_f$ & 11000 & $s$ & final time for study 1\\
\hline
$c_\alpha^+$ &0.521& $kg/kg$ water & fraction of $\alpha$-lactose in feed\\
\hline $c_\beta^+$ &0.359& $kg/kg$ water& fraction of $\beta$-lactose in feed\\
\hline
$V_0$ & 0.0015 & $m^3$  &initial volume\\
\hline
$V_{\rm max}$ & 0.01 & $m^3$  &maximum  volume\\
\hline\hline
\end{tabular}
\caption{Numerical constants}
\end{center}
\end{table}

\begin{table}[h!]
\begin{center}
\begin{tabular}{||c|c|c||}
\hline \hline
$n(L,t)$ &  $\sharp/m.m^3$  &particle size distribution\\
\hline
$m_\alpha(t)$& $kg$ & mass of $\alpha$-lactose in solution\\
\hline
$m_\beta(t)$ &$kg$& mass of $\beta$-lactose in solution\\
\hline
$V(t)$ & $kg$ & volume of slurry\\
\hline
$A(t)$ & $m^2$ & contact surface \\
\hline\hline
\end{tabular}
\caption{Units of dynamic quantities}
\end{center}
\end{table}

\section{Optimal control problem}
\label{optimal-problem}
The benefit of the moment approach 
is that we may choose a finite number of moment equation to replace
(\ref{pop}). Our present approach is to retain a sufficient number of moments so that
the salient features of the seed $n_0(L)$ 
may be captured by these moments, and in our experiments
we decided to retain the moments $\mu_0,\dots,\mu_5$. The remaining states
of the system dynamics are then $m_{\rm H_2O}$, $m_{\rm cry}$, $c_\alpha$,
$c_\beta$, $T$, $T_{\rm jacket}$. The control
inputs are $u_1=T_{\rm sp}$ and $u_2=q_{\rm H_2O}$.

\begin{table}[h!]
\begin{center}
\begin{tabular}{||c|c||c|c||}
\hline\hline
$\mu_0$ &$1.2405 1^{10}$   & $m_{\rm H_2O}$& $0.92$kg\\
\hline
$\mu_1$ & $2.1767 10^{6}$  &$c_\alpha$ &0.359\\
\hline
$\mu_2$ & $409.2491$ &$c_\beta$ & 0.521\\
\hline
$\mu_3$ & $0.0812$   &$T$ &$70^0$C \\
\hline
$\mu_4$ & $1.6812 10^{-5}$  &$T_{\rm jacket}$ &$20^0$C \\
\hline
$\mu_5$ & $3.6094 10^{-9}$ & &\\
\hline\hline
\end{tabular}
\end{center}
\caption{Initial values for study 1}
\end{table}

In this work, we compare between several policies :
\begin{itemize}
\item Policy 1 : We fix the values of the set-point temperature 
 $T_{\rm sp}=15  [{\,}^0C]$ and the feed rate of solvent $q_{\rm H_2O}=0.0056  [Kg/h]$.  This policy is
 referred to as {\em constant} in the figures.
\item Policy 2 : Here we fix  the value of the feed rate of solvent 
  $q_{\rm H_2O}=0.0056  [Kg/h]$,  while the  set-point temperature $T_{\rm sp}(t)$ starts
  at $T_{\rm sp}(0)=15  [{\,}^0C]$ and decreases
linearly. This policy is called {\em linear} in the figures.
\item Policy 3 :   We  control  the set-point temperature 
 $u_1(t)=T_{\rm sp}(t)$ and also the  feed rate of solvent $u_2(t)=q_{\rm H_2O}(t)$ using various objectives. 
  This policy is called {\em optimal} in the figures.
\end{itemize}

\subsection{Scenario 1}
\label{Scenario-1}
Our first control problem minimizes the weighted mean size diameter $d_{43}=\dfrac{\mu_{4}}{\mu_{3}}$
at  fixed final time $t_f=11 000$ seconds.
This is cast as the optimization program
\begin{eqnarray}
\label{scenario1}
\begin{array}{ll}
\mbox{minimize} & d_{43}(t_f)=\displaystyle\frac{\mu_4(t_f)}{\mu_3(t_f)} \\
\mbox{subject to} 
&\mbox{dynamics }(\ref{mom1}) - (\ref{dmcry}) \\
& 0 \leq V(t) \leq V_{\max} \\
& 0^0{\rm C} \leq T(t) \leq 70^0{\rm C} \\
&c_\alpha(t) \ge c_{\alpha,\rm sat}(c_\beta(t),T(t)) \\
&0^0{\rm C} \leq T_{\rm sp}(t) \leq 40^0{\rm C} \\
&0 \leq q_{\rm H_2O}(t) \leq 0.1
\end{array}
\end{eqnarray}
The control variables are set-point temperature $u_1(t)=T_{\rm sp}(t)$ and feed rate of water $u_2(t)=q_{\rm H_2O}(t)$.
The percentages $\dot{c}_\alpha^+$ and $\dot{c}_\beta^+$ of lactose in the feed are
kept constant.

\subsection{Scenario 2}
\label{Scenario-2}
Our second control problem minimizes the nucleation rate $B(c_\alpha,c_\beta,T)$
at the fixed final time $t_f=11 000$ seconds.
This is cast as the optimization program

\begin{eqnarray}
\label{scenario2}
\begin{array}{ll}
\mbox{minimize} & B(t_f) \\
\mbox{subject to} & \mbox{constraints of }(\ref{scenario1}) \\
\end{array}
\end{eqnarray}
The control variables are again $T_{\rm sp}$ and $q_{\rm H_2O}$.


\subsection{Scenario 3}
\label{Scenario-3}
Our third control problem minimizes the coefficient of variation $CV$
at the fixed final time $t_f=11 000$ seconds.
This is  the optimization program
\begin{eqnarray}
\label{scenario3}
\begin{array}{ll}
\mbox{minimize} & CV(t_f) = \dfrac{\mu_3(t_f)\mu_5(t_f)}{(\mu_4(t_f))^{2}} - 1 \\
\mbox{subject to} & \mbox{constraints of }(\ref{scenario1}) 
\end{array}
\end{eqnarray}
The control variables are again $T_{\rm sp}$ and $q_{\rm H_2O}$.

\begin{figure}[!ht]
\centering
\hspace*{-0.7cm}
\includegraphics[width=0.58\textwidth]{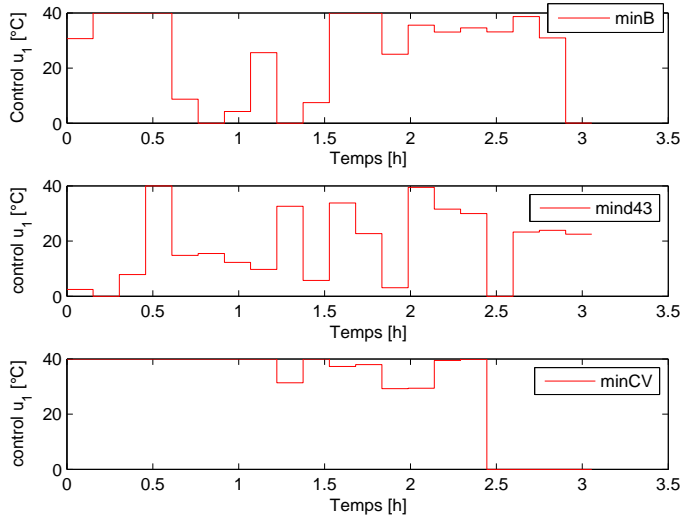} 
 \caption{Optimal set-point temperature profile $u_1=T_{\rm sp}$  for  
 minimization of the three criteria $B$, $CV$ and $d_{43}$ with fixed final time $t_f$.}
 \label{Scen_temp}
\end{figure}

 \begin{figure}[!ht]
 \centering
\hspace*{-0.5cm}
\includegraphics[width=0.58\textwidth]{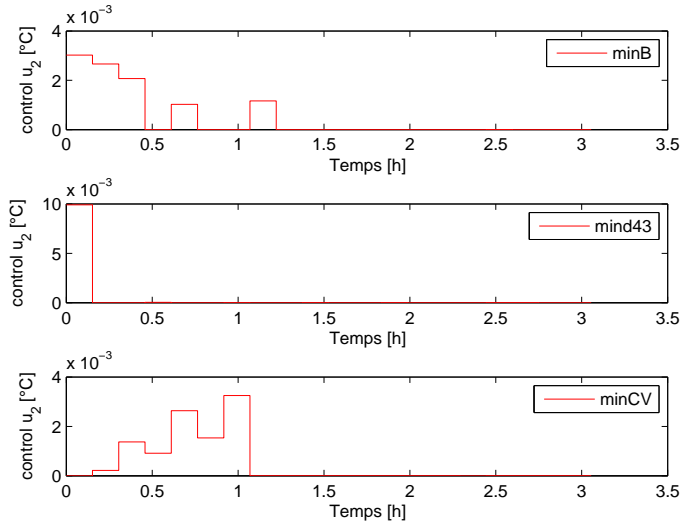}
\caption{Optimal feed profile $u_2 = q_{\rm H_2O}$  for minimization of
the three criteria $B$, $d_{43}$ and $CV$ with fixed final time $t_f$. }
 \label{Scen_feed}
 \end{figure}

 \begin{figure}[!ht]
\hspace*{-0.5cm}
\includegraphics[width=0.5\textwidth]{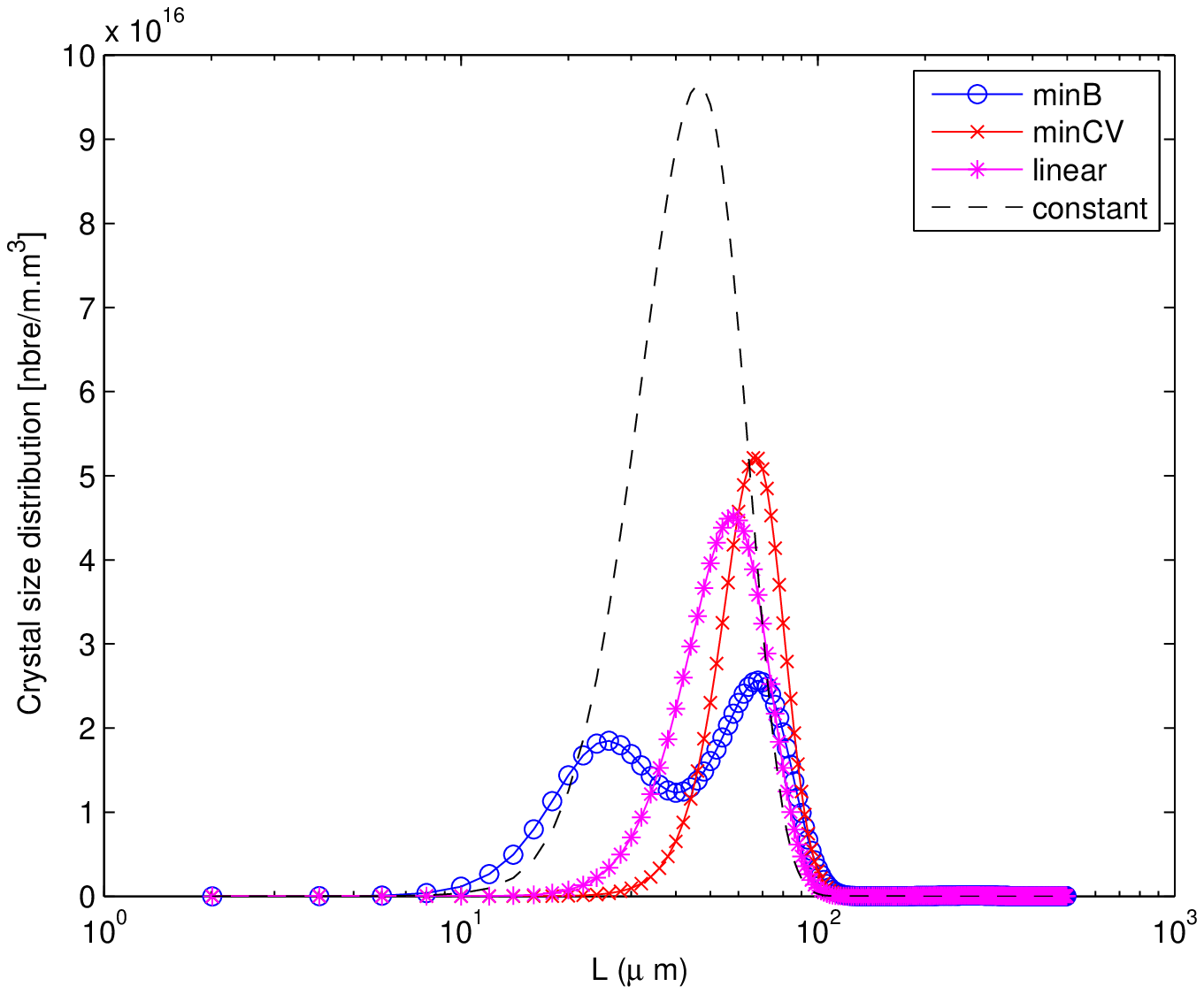}
\includegraphics[width=0.5\textwidth]{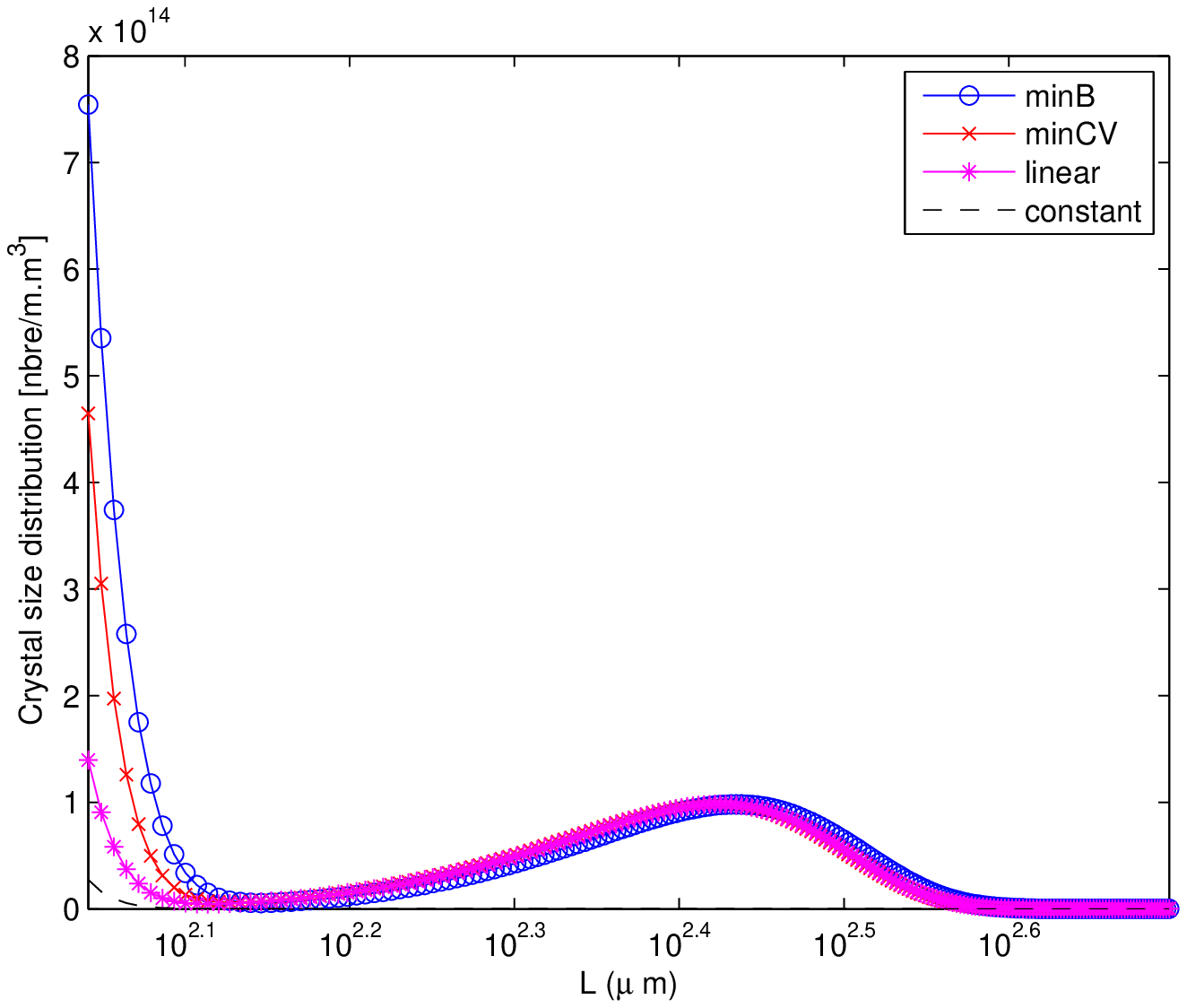}
\caption{Final crystal size distribution $L \mapsto n(L,t_f)$ displayed for minimization of $B$ (blue) and $CV$ (red)
compared with linear (magenta)  and constant (dashed black) policies for fixed final time  $t_{\text{final}}$. Right hand image
shows zoom on range $[10^{2.1},10^{2.6}]$.}
 \label{Scen_dist1}
\end{figure}

 \begin{figure}[!ht]
\hspace*{-0.5cm}
\includegraphics[width=0.5\textwidth]{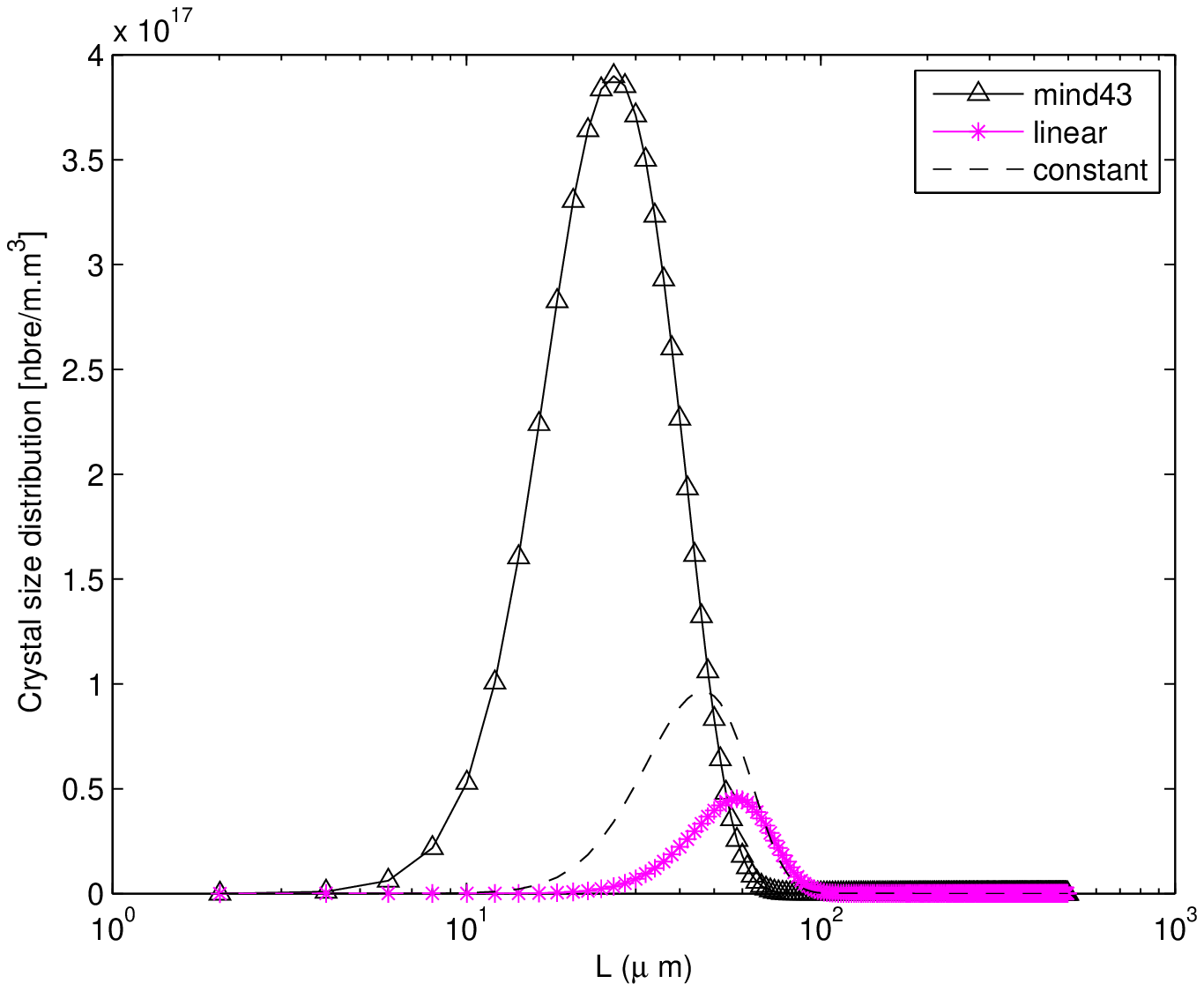}
\includegraphics[width=0.5\textwidth]{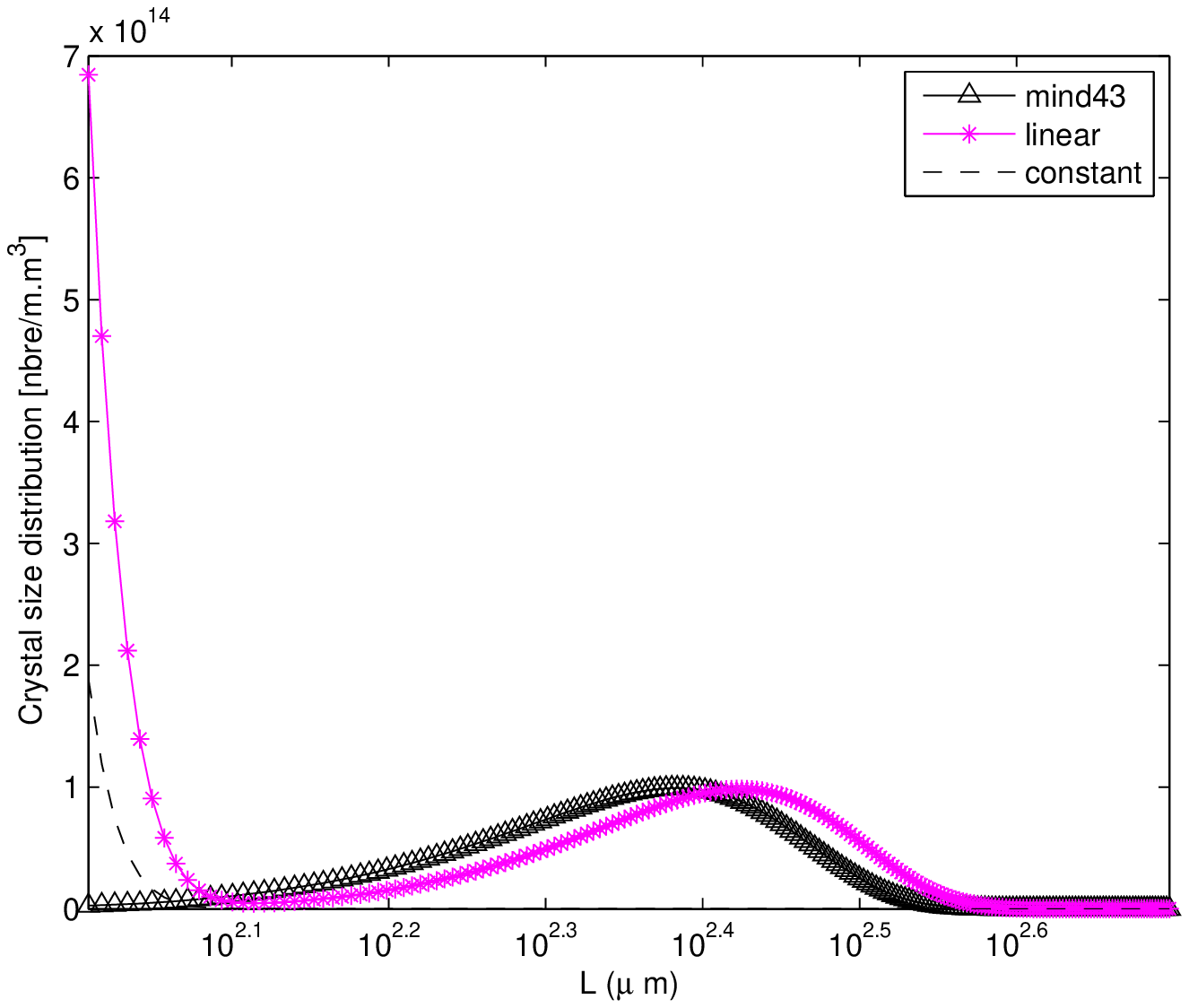}
\caption{Final crystal size distribution $L \mapsto n(L,t_f)$ for minimizing $d_{43}$ (black solid), 
compared with linear (magenta) and constant (black dashed) policies for fixed final time  $t_{\text{final}}$. Right hand
image shows zoom on range $[10^{2.1},10^{2.6}]$.}
 \label{Scen_dist2}
\end{figure}

 \begin{figure}[!ht]
 \hspace{-0.4cm}
 \includegraphics[width=0.5\textwidth]{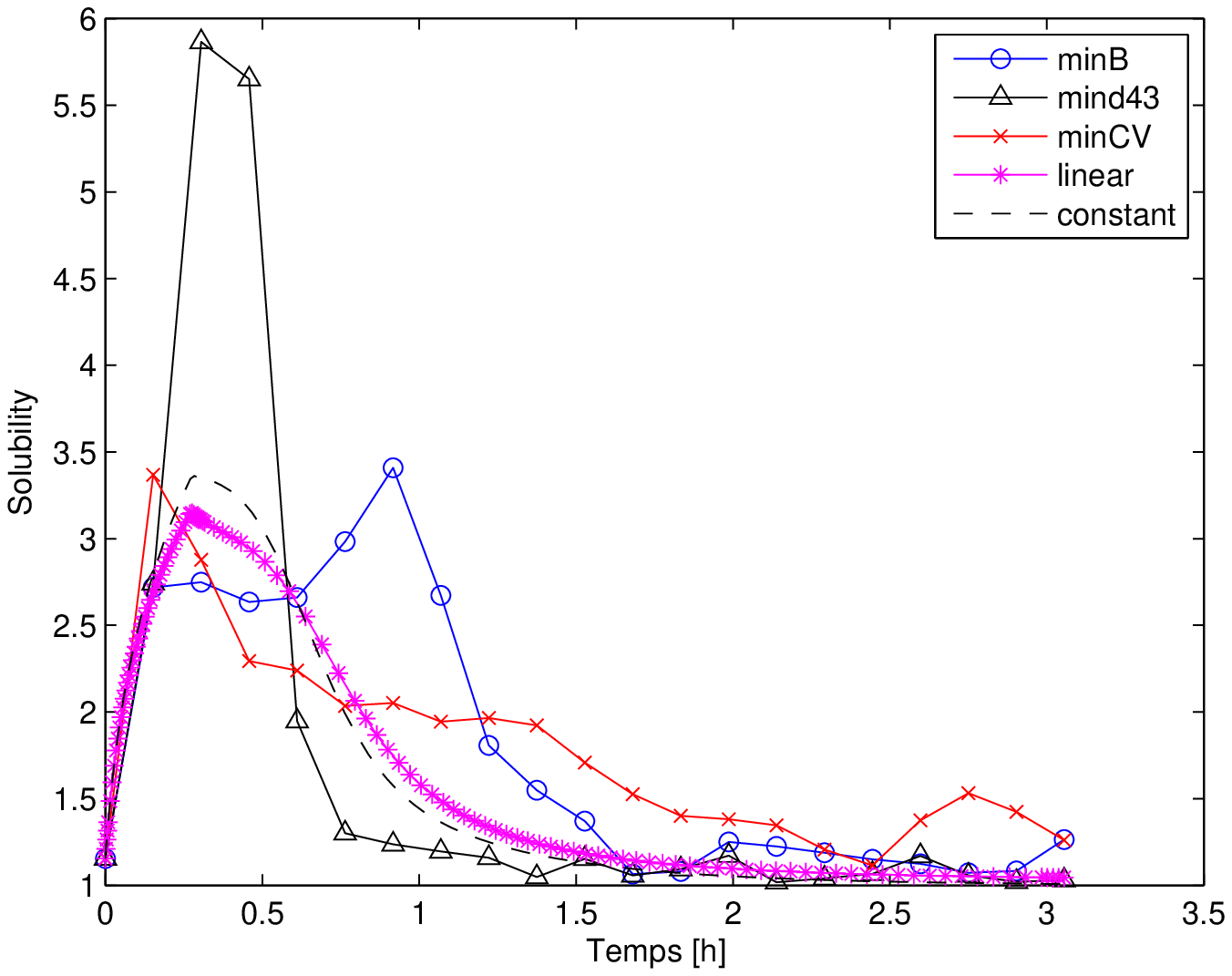}
 \includegraphics[width=0.5\textwidth]{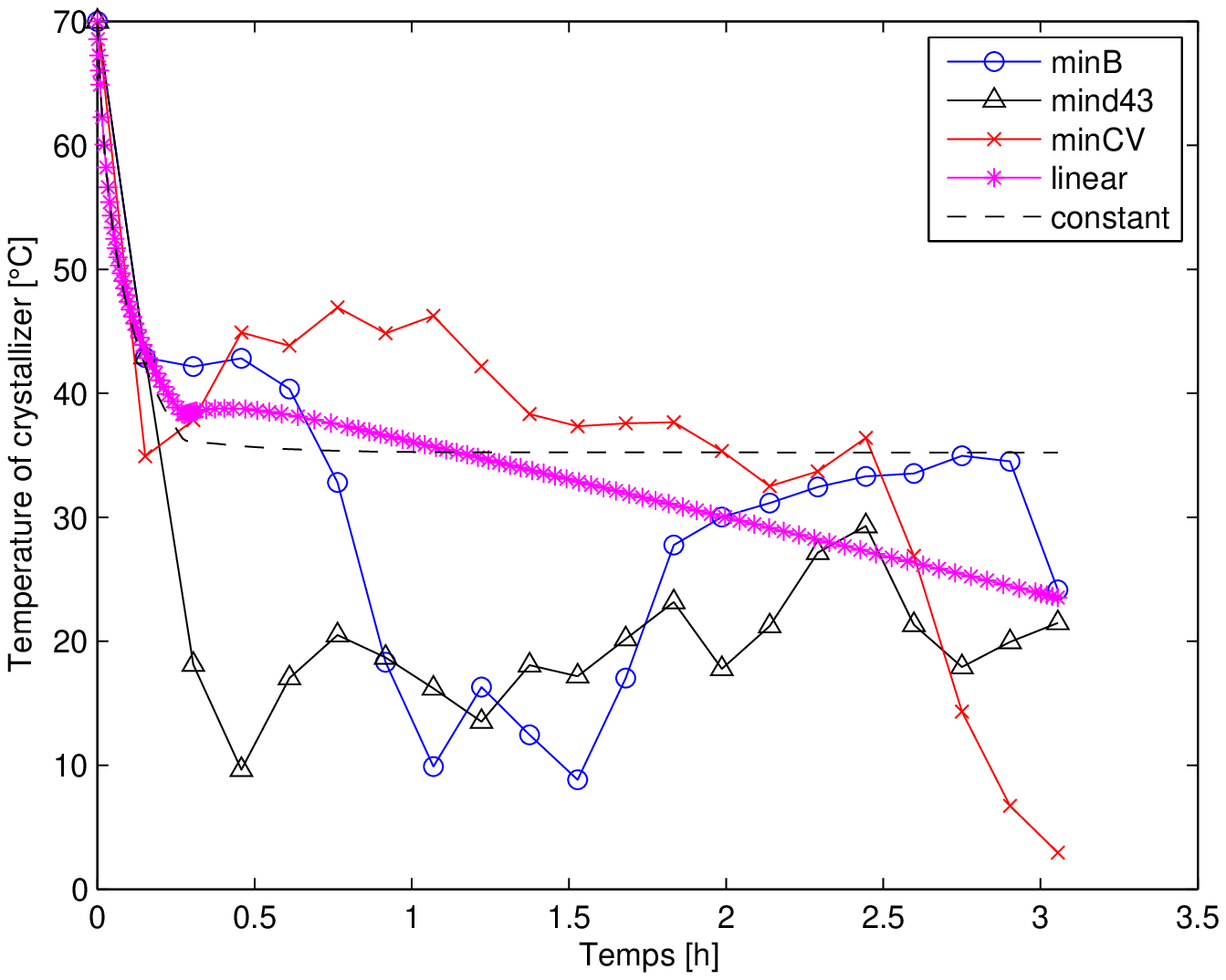}\\
\caption{Left image compares solubility,  right image compares temperature of crystallizer  for the five
policies constant (dashed black), linear (magenta), optimal with $B$ (blue), optimal with $CV$ (red) and 
optimal with $d_{43}$
(black continuous)  for  $t_{\text{final}}$ fixed. }
 \label{Scen_solub}
\end{figure}

 \begin{figure}[!ht]
 \hspace{-0.4cm}
 \includegraphics[width=0.5\textwidth]{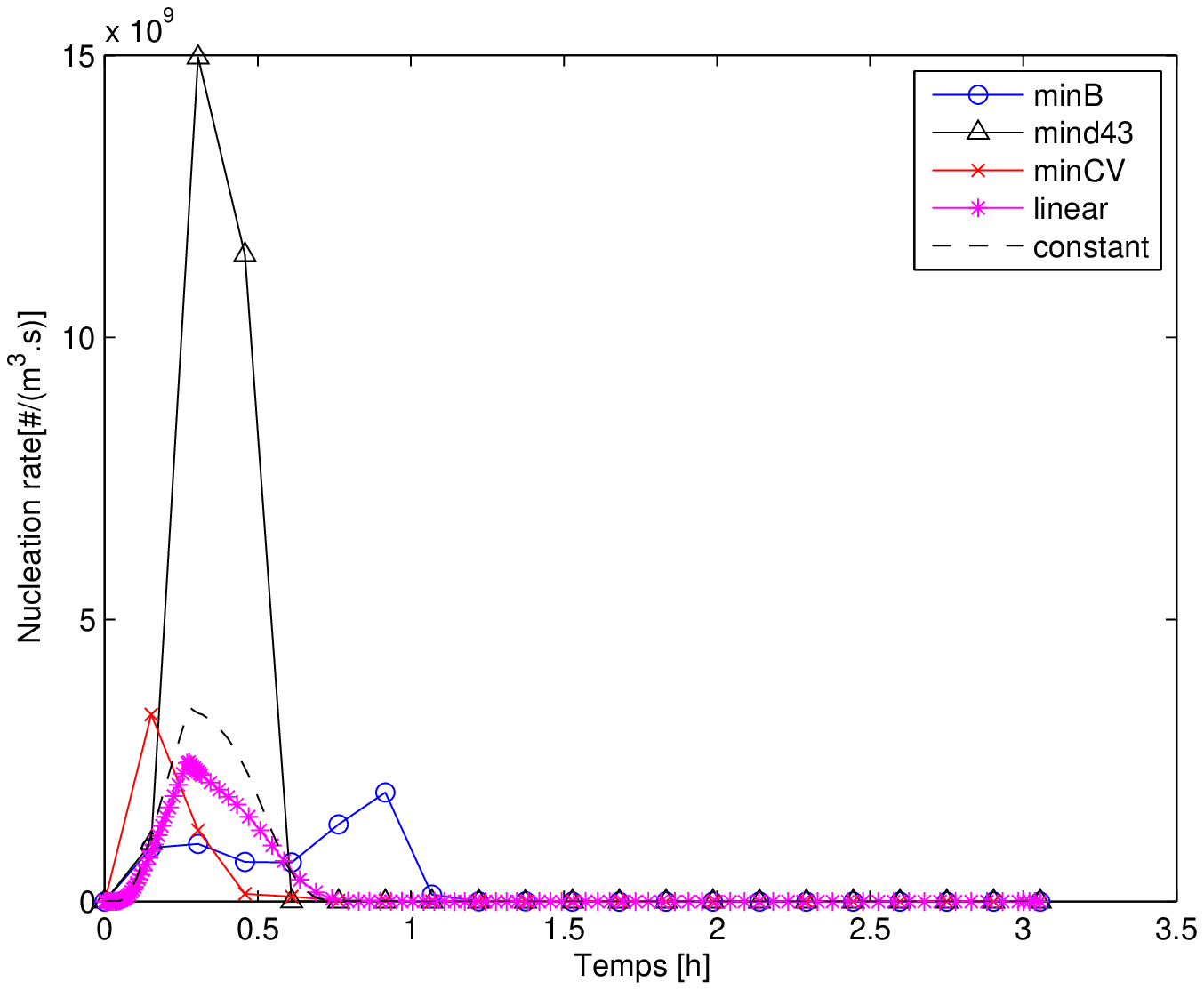}
 \includegraphics[width=0.5\textwidth]{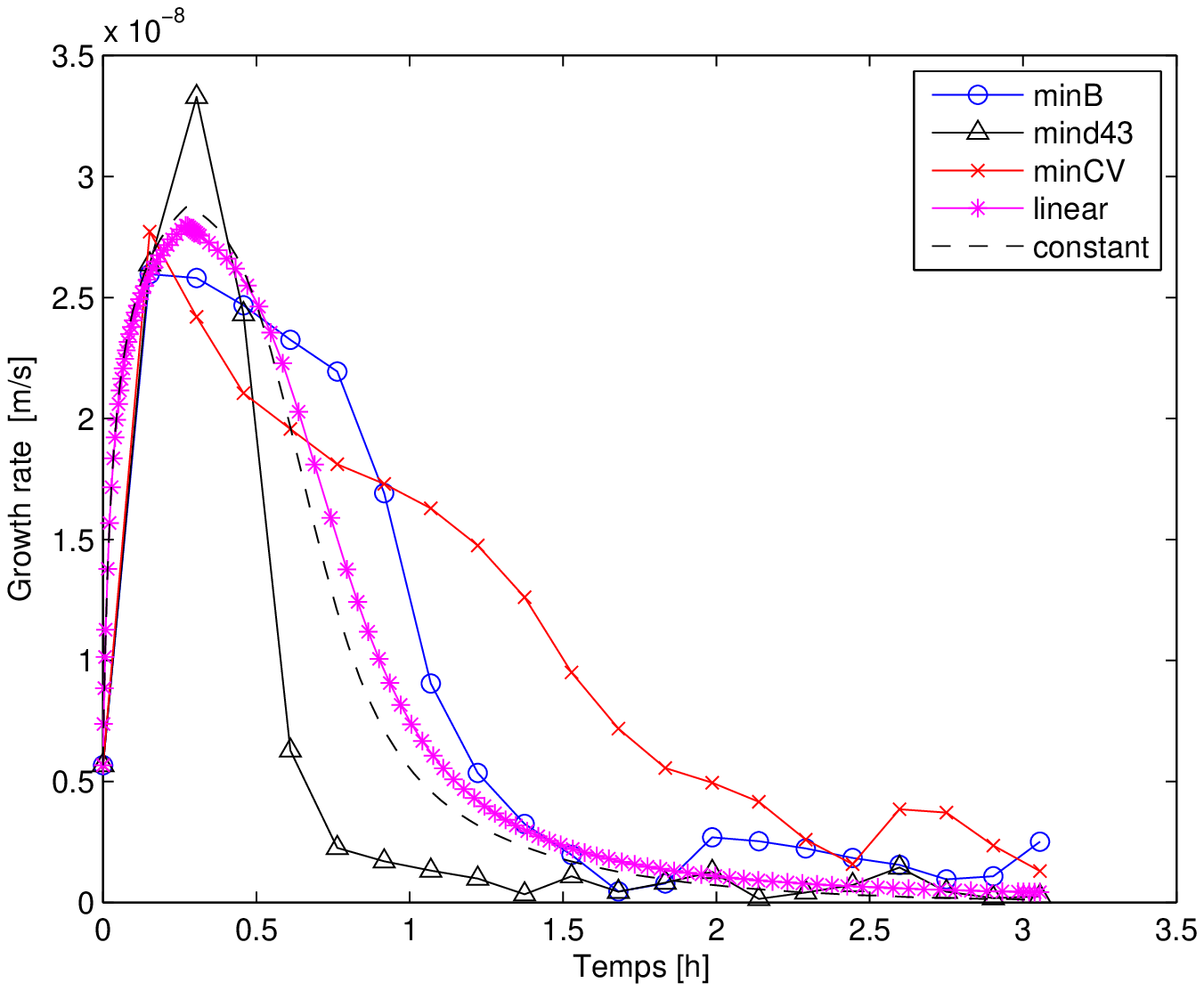}
\caption{Left image compares  growth rate, right image compares  crystal mass for the five policies
constant (dashed black), linear (magenta), optimal with $B$ (blue), optimal with $CV$ (red),
optimal with $d_{43}$ (black continuous) for  $t_{\text{final}}$ fixed. }
 \label{Scen_grow}
\end{figure}

\begin{figure}[!ht]
 \hspace{-0.4cm}
 \includegraphics[width=0.5\textwidth]{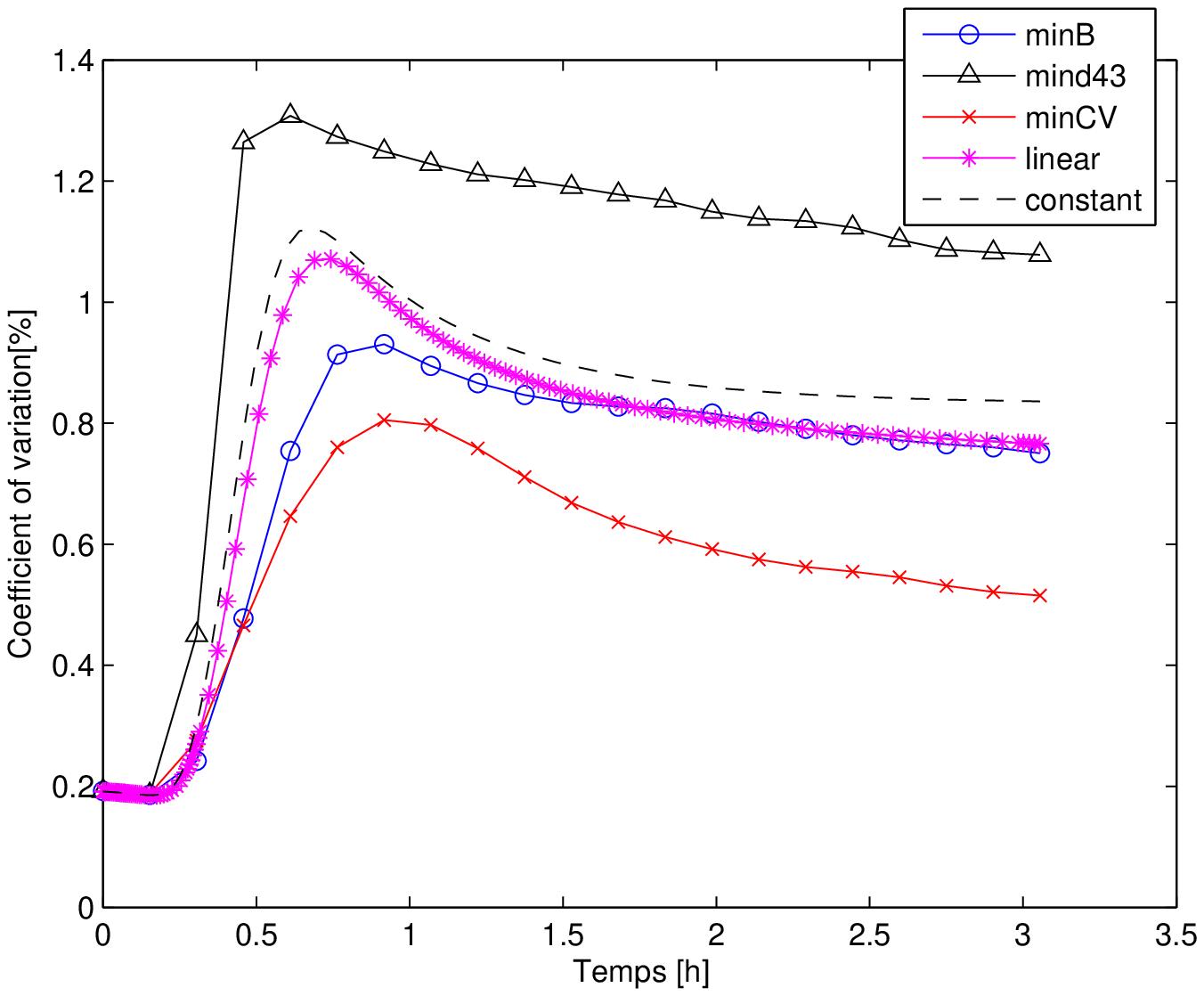}
 \includegraphics[width=0.5\textwidth]{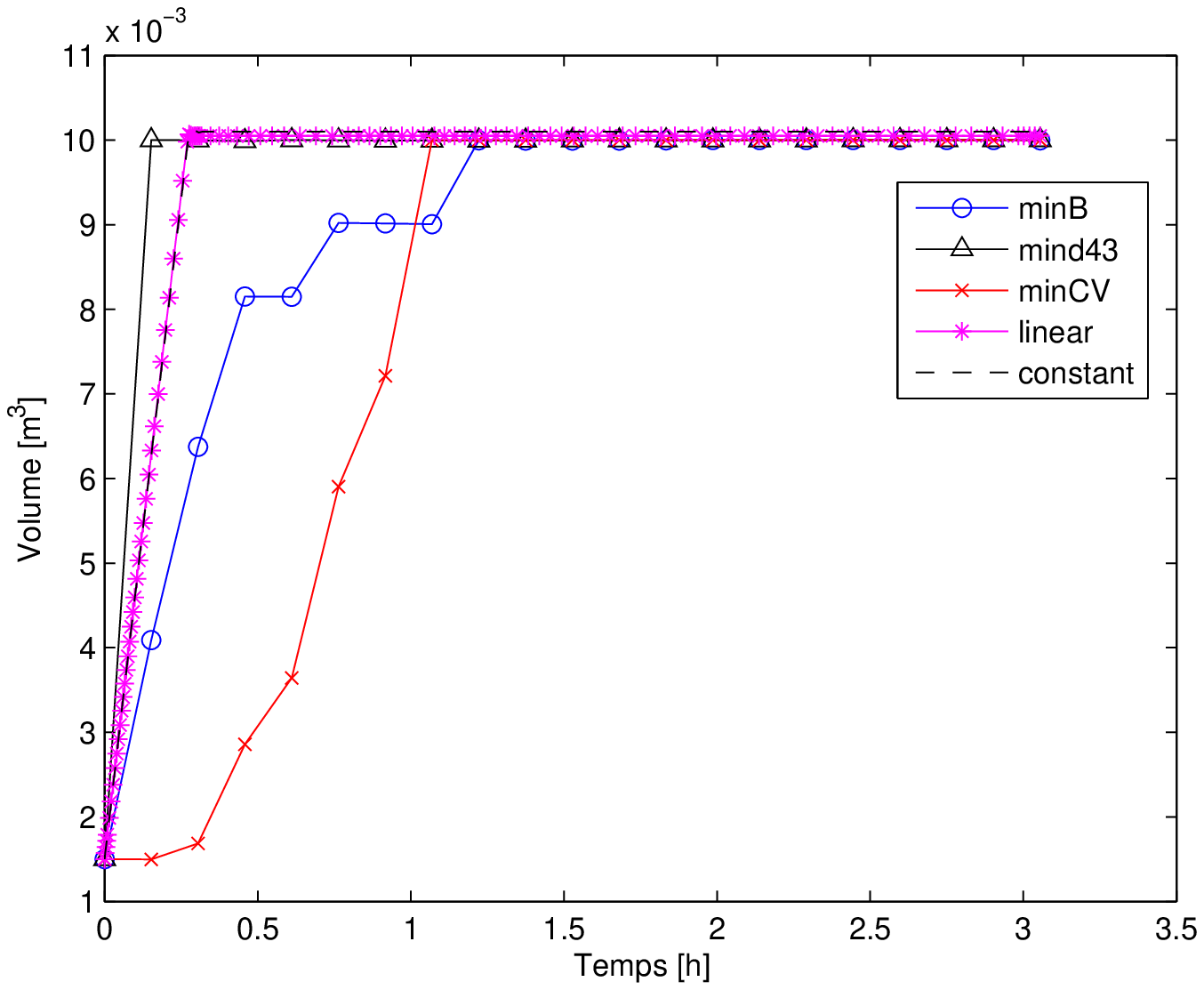}
\caption{Evolution of   coefficient of variation $CV(t)$ left, evolution of volume $V(t)$   right.  Comparison
 of the five policies constant (dashed black), linear (magenta), optimal for $B$ (blue),
 optimal for $CV$ (red), optimal for $d_{43}$ (black continuous), for  $t_{\text{final}}$ fixed. }
 \label{Scen_cv}
\end{figure}

\begin{figure}[!ht]
 \hspace{-0.4cm}
 \includegraphics[width=0.5\textwidth]{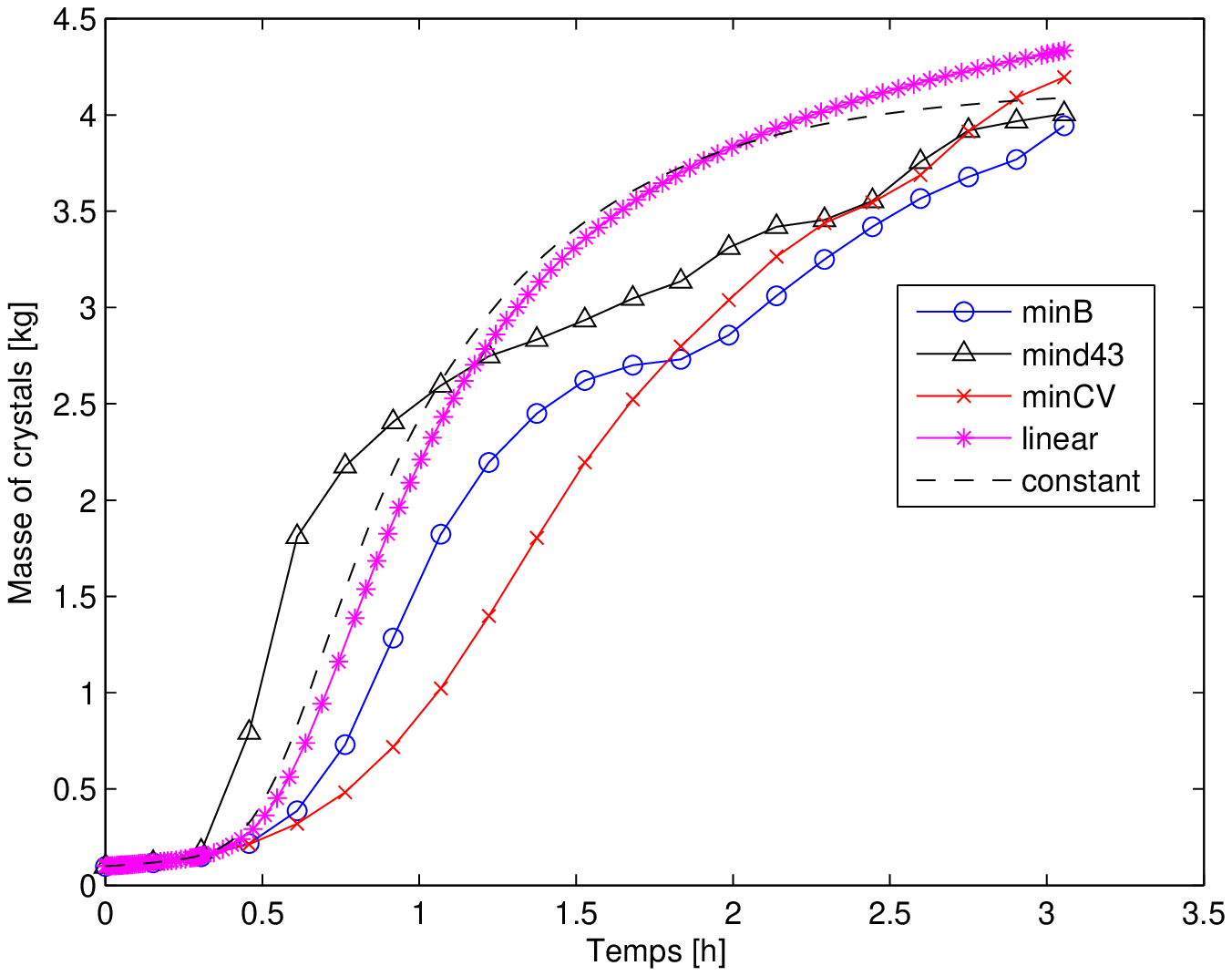}
 \includegraphics[width=0.5\textwidth]{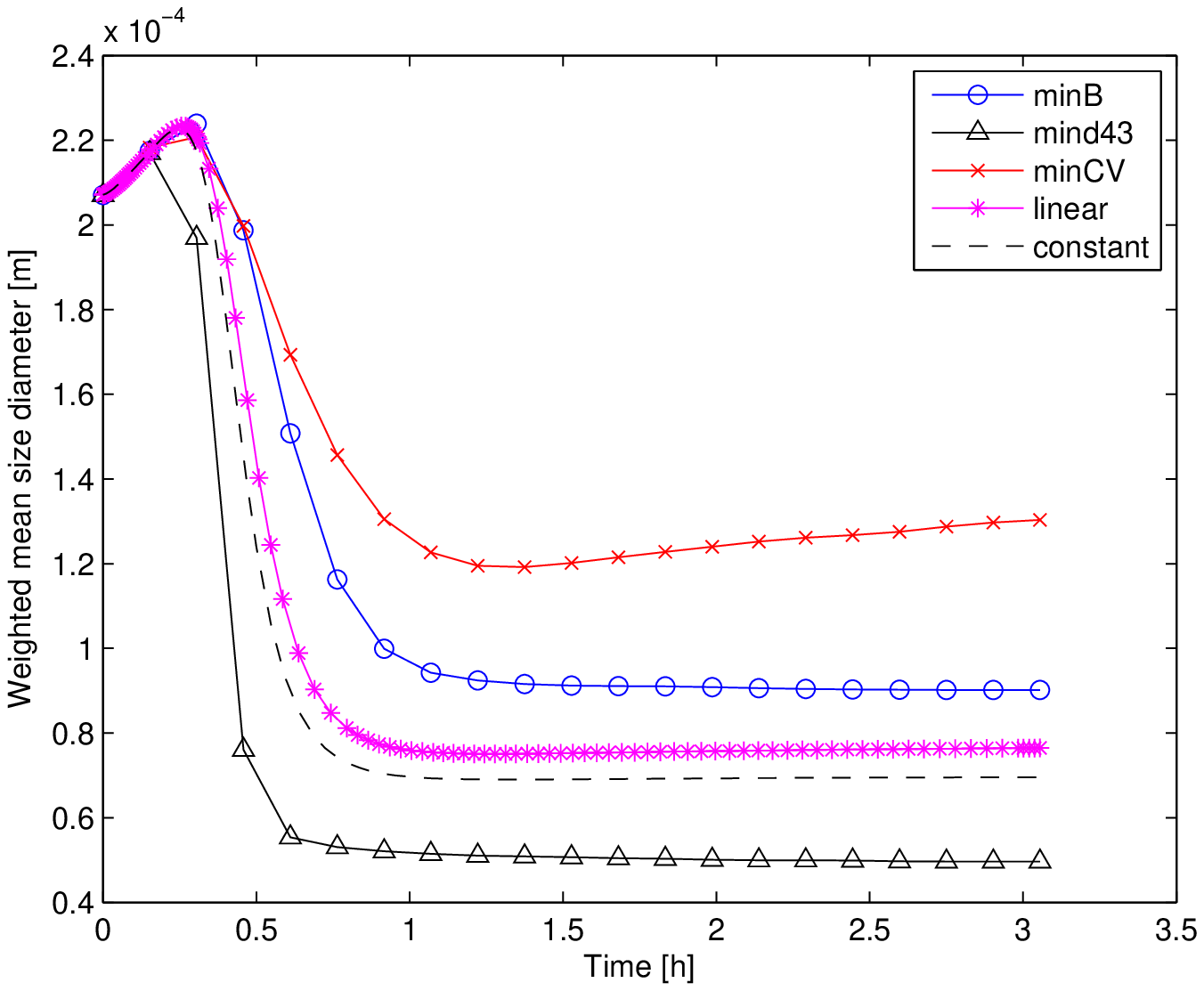}
\caption{Evolution of  weighted mean size diameter  $d_{43}$ left, evolution of crystal mass   right.  Comparison
 of the five policies constant (dashed black), linear (magenta), optimal for $B$ (blue),
 optimal for $CV$ (red), optimal for $d_{43}$ (black continuous), for  $t_{\text{final}}$ fixed. }
 \label{Scen_d43}
\end{figure}

In Figures \ref{Scen_temp}, \ref{Scen_feed}  we present results obtained with the optimal regulation 
of set-point temperature $u_1=T_{\rm sp}$ and  feed
rate $u_2=q_{\rm H_2O}$ and compare these  to more standard scenarios, where
 temperature and feed rate are fixed or follow simple heuristic profiles proposed in the literature.

In  Figure \ref{Scen_dist1} we present the optimal
crystal size distribution obtained from minimization of the  weighted mean size diameter  $d_{43}$
in (\ref{scenario1}) 
compared with standard scenarios.  
Figure \ref{Scen_dist2} shows  the crystal size distribution for the optimal control of  nucleation rate  $B$ 
and coefficient of variation $CV$ compared to the more standard scenarios. The optimal profile for the nucleation 
rate shows the existence of 
two peaks which indicates the existence of two populations of crystals.

In  Figure \ref{Scen_solub}  we present the evolution of solubility of $\alpha-$lactose  and the temperature of 
the crystallizer 
by comparing
several scenarios. The profile of solubility shows an early peak, which correspond 
to a sharp decrease in the  temperature profile of the crystallizer.  
Comparison between the cost functions shows that the highest peak occurs when minimizing 
the weighted mean size diameter  $d_{43}$. 
In the case of minimization of nucleation rate $B$, we see the existence of two peaks which correspond with two peaks on crystal size distribution profile.

In  Figure \ref{Scen_grow} we present the evolution of nucleation rate $B$ and  growth rate $G$ in 
comparison between the several scenarios. At the beginning of the profile of nucleation rate $B$, we note the highest peak in case of minimization of  weighted mean size diameter  $d_{43}$.

In  Figure \ref{Scen_cv}  we present the coefficient of variation $CV$ and  volume of crystallizer  in 
comparison between all objectives and scenarios. The volume profiles show that optimization
of different cost functions may lead to fairly different ways of 
filling the crystallizer in the semi-batch phase. For instance, filling in the  linear policy occurs much faster 
than  for  minimization of the coefficient of variation, which gives the slowest filling.

Figure \ref{Scen_d43} compares the profiles of overall crystals mass and of the 
weighted mean size diameter  $d_{43}$ in all scenarios.

\subsection{Scenario 4}
The next extension is to add the moments
of $n_0(L)$ as unknown parameters, the idea being that a suitable choice
of the initial seed of given mass should give even better results.
We decide to fix the total mass of crystal seed as $k_v V_0 \rho_{\rm cry}\int_0^\infty n_0(L)L^3 dL = 0.1$kg.
That leads to the optimization program
\begin{eqnarray}
\begin{array}{ll}
\mbox{minimize} & d_{43} \\
\mbox{subject to} & \mbox{constraints of } (\ref{scenario2}) \\
&\mu_3(0) =\displaystyle \int_0^\infty n_0(L)L^3dL =0.0812
\end{array}
\end{eqnarray}
where now $T_{\rm ref}(t)$, $q_{\rm H_2O}(t)$ and $\mu_0(0)$, $\mu_1(0)$,
$\mu_2(0)$, $\mu_4(0)$, $\mu_5(0)$ are optimization variables. 

The interesting point of this program is that
once the optimal solution $(T_{\rm ref}^*,q_{\rm H_2O}^*,\mu_\nu^*)$ is reached, we need to reconstruct
a function $n_0^*(L)$ such that  its moments $0,\dots,5$ coincide with
$\mu_0^*,\mu_1^*,\mu_2^*,\mu_3^*=0.0812, \mu_4^*,\mu_5^*$. This can be achieved by solving
the maximum entropy function reconstruction problem
\begin{eqnarray}
\label{maxent}
\begin{array}{ll}
\mbox{minimize} & \displaystyle \int_0^\infty n_0(L) \log n_0(L) dL \\
\mbox{subject to} &\displaystyle \int_0^\infty L^\nu n_0(L) dL = \mu_\nu^*, \nu = 0,\dots,5.
\end{array}
\end{eqnarray}
Notice that (\ref{maxent}) may be solved
by standard software, see e.g. Borwein {\em et al.} \cite{browein1}, \cite{browein2}.

\subsection{Scenario 5}
The natural figure of merit to maximize the crystal mass
within a certain range $L_1 \leq L \leq L_2$ of small particle sizes is
\begin{eqnarray}
\label{natural}
\max \int_{L_1}^{L_2} L^3 n(L,t_f)dL
\end{eqnarray}
at the final time $t_f$, but this objective is not directly
accessible in the moment approach. Substrates like $B,CV,d_{43}$ are
non-specific and
must be expected to give only a  crude approximation of (\ref{natural}). We therefore
propose
the following more sophisticated strategy, which is compatible with the moment approach.

We define a target particle size distribution $n_1(L)$, which has a bulk in the range $[L_1,L_2]$,
normalized to satisfy
\[
\nu_3 = \int_0^\infty L^3 n_1(L) dL = 1.
\] 
Then  we compute as many of its moments $\nu_0,\dots,\nu_N$
as we wish to use, where as before $N=5$ in our tests. The optimization
program we now solve is
\begin{eqnarray}
\label{scenario4}
\begin{array}{ll}
\mbox{minimize} & \displaystyle\sum_{i=0}^N w_i \left( \mu_i(t_f) - \mu_3(t_f) \nu_i\right)^2 \\
\mbox{subject to}& \mbox{constraints of } (\ref{scenario1}) \\
&t_f \leq t_{\rm max}
\end{array}
\end{eqnarray}
where the $w_i$ are suitably chosen weights. Notice that the least squares objective of (\ref{scenario4})
tries to bring the moments of $n(L,t_f)$ as close as possible to the
moments of the function $\mu_3(t_f)n_1(L)$, which has the correct shape, and has the same
total crystal  mass as $n(L,t_f)$. {Here the final time is considered free}.

\begin{figure}[h!]
\hspace{-.4cm}
\includegraphics[width=0.5\textwidth]{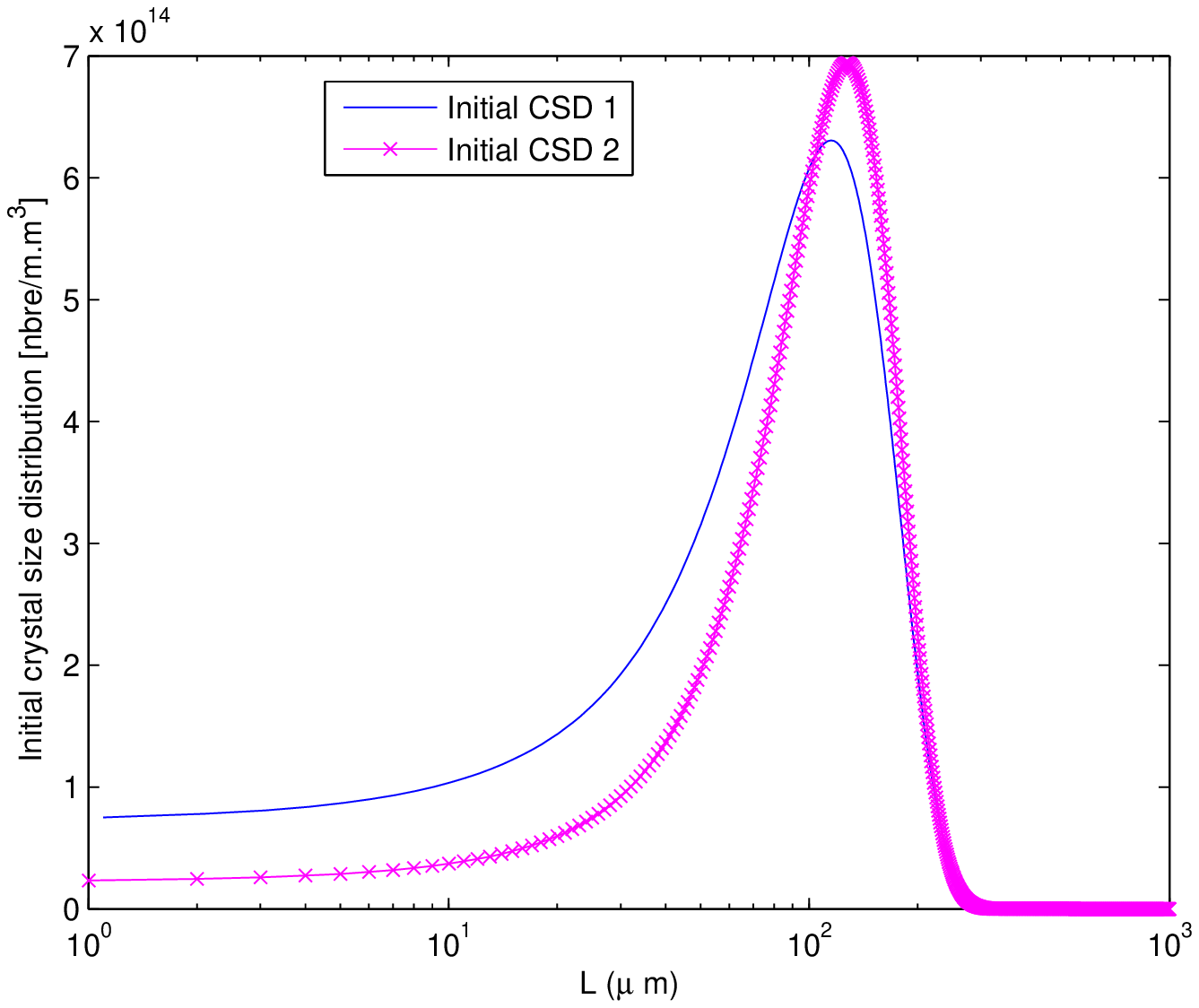}
\includegraphics[width=0.5\textwidth]{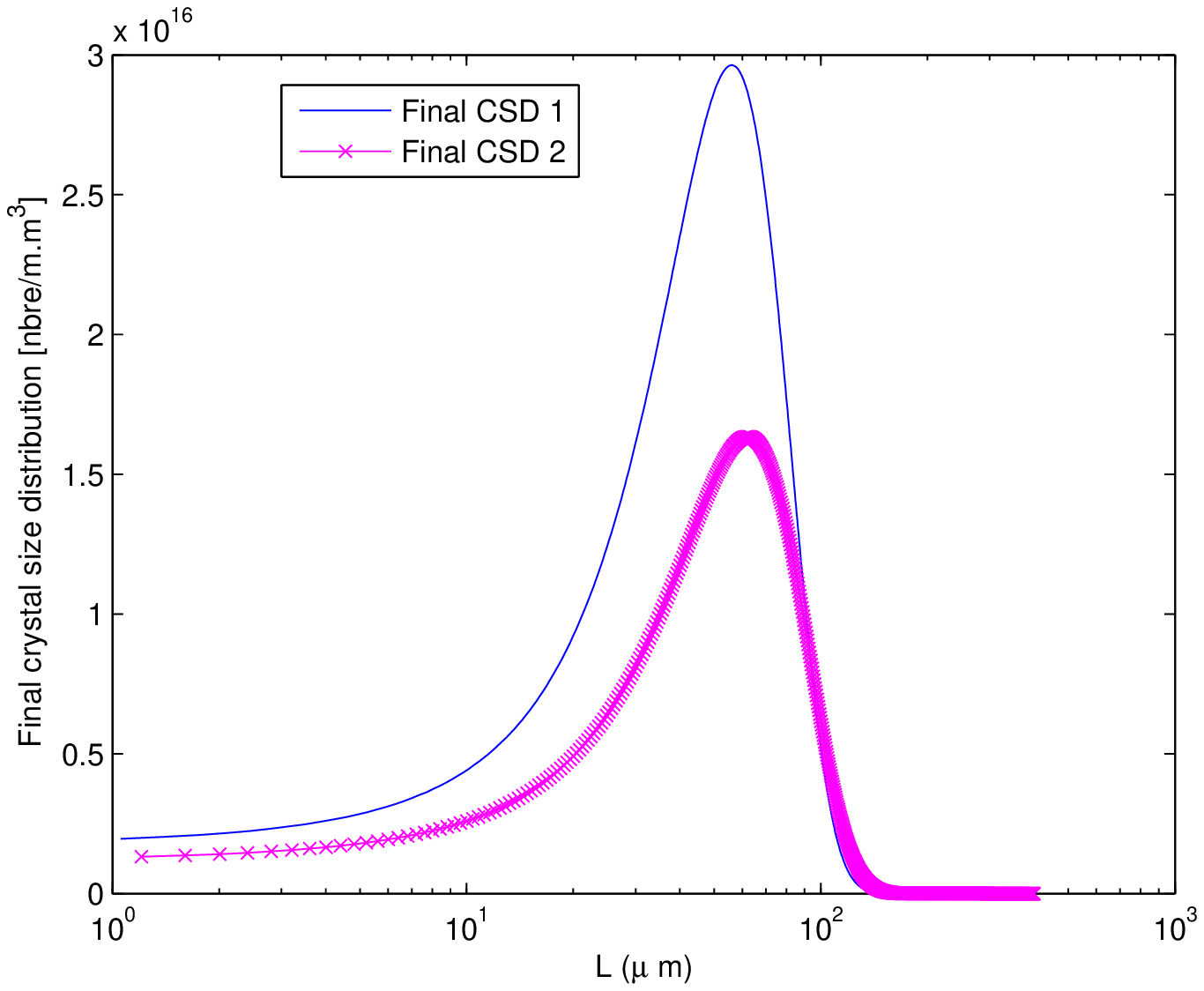}
\caption{Scenario 5: Different initial seeds $n_0(L)$ shown on the left lead to different products $n(L,t_f)$ on the right.
Blue uses $\mu_3(0)=0.291$, magenta uses $\mu_3(0)=0.401$.}
\end{figure}

\section{Method}
For our testing we have used the solver ACADO \cite{acado} based on a semi-direct 
single or multiple-shooting strategy, including automatic differentiation, based ultimately
on the semi-direct multiple-shooting algorithm of Bock and Pitt \cite{bock}. ACADO is a 
self-contained public domain software environment 
written in {\tt C++} for automatic control and dynamic optimization. 

Alternatively, we also use the solver PSOPT \cite{psopt}, which is a public domain  extension of the
NLP-slover IPOPT \cite{ipopt} or SNOPT \cite{snopt}  and is based on pseudospectral optimization 
which uses Legendre or Cheybyshev polynomials and discretization  based on Gauss-Lobatto nodes.

A difficulty with both solvers is the strong dependence of convergence and solutions on the initial guess,
as must be expected in a local optimization context. Our testing shows that it is often mandatory
to have a simulated study $(x_{\rm init},u_{\rm init})$ available to start the optimization
from that point. This initial guess may use parameters  from a previous optimization
study, which give already a decent cost in the present study. In some cases homotopy techniques, using for instance $t_f$ as a parameter,
have to be used.

Once optimal controls $u^*=(T_{\rm ref}^*,q_{\rm H_2O}^*)$ have been computed in
any  one of the scenarios, 
we use the full crystallizer model
(\ref{pop}), (\ref{bdry}), (\ref{water}) --  (\ref{eq4})
to simulate the system, using an initial seed $n_0(L)$
which produces the initial moments $\mu_\nu(0)$. 
In those cases where the moments of the initial seed are parameters, which are also optimized,
we use the optimal $\mu^*=(\mu_0^*,\dots,\mu_N^*)$ to compute an estimation
$n_0^*(L)$ of the optimal crystal seed with these moments using \cite{browein1} and \cite{browein2}.

The final stage in each experiment is a simulation of the full population balance
model using the optimal $(T_{\rm sp}^*,q_{\rm H_2O}^*,n_0(L)^*)$, obtained from the
moment-based optimal control problem.

\section{Conclusion}
We have presented and tested several control strategies which allow
to maximize the crystal mass of particles of small size, typically in a range
of $10^{-5}-10^{-4} \mu m$. Our approach
was compared to more standard heuristic control policies used in the literature to  
regulate  temperature and feed rate in semi-batch mode. 
Our simulated numerical results show that  it is beneficial to apply 
optimal control strategies in semi-batch solvated crystallization of $\alpha$-lactose monohydrate.

\end{document}